\documentclass[a4paper,11pt] {amsart}
\usepackage{amssymb,latexsym,amsmath}
\usepackage{amsfonts,amsbsy,bm}
\def\a {\alpha }

\def\O{\Omega}

\def\a{\alpha}

\def\Da{\Delta}
\def\ga{\gamma}

\def\ba{\beta}

\newcommand{\bprop} {\begin{proposition}}
\newcommand{\eprop} {\end{proposition}}
\newcommand{\btheo} {\begin{theorem}}
\newcommand{\etheo} {\end{theorem}}
\newcommand{\blem} {\begin{lemma}}
\newcommand{\elem} {\end{lemma}}
\newcommand{\bcor} {\begin{corollary}}
\newcommand{\ecor} {\end{corollary}}

\newcommand{\Be}{\begin{equation}}
\newcommand{\Ee}{\end{equation}}
\newcommand{\Bea}{\begin{eqnarray}}
\newcommand{\Eea}{\end{eqnarray}}
\newcommand{\Bes}{\begin{equation*}}
\newcommand{\Ees}{\end{equation*}}
\newcommand{\Beas}{\begin{eqnarray*}}
\newcommand{\Eeas}{\end{eqnarray*}}
\newcommand{\Ba}{\begin{array}}
\newcommand{\Ea}{\end{array}}

\def\C{\mathbb{C}}

\newtheorem{theorem}{Theorem}[section]
\newtheorem{corollary}[theorem]{Corollary}
\newtheorem{lemma}[theorem]{Lemma}
\newtheorem{proposition}[theorem]{Proposition}

\newtheorem{remark}[theorem]{Remark}

\begin{document}
\title[Off-diagonal Estimates]{Off-diagonal estimates of some Bergman-type operators on tube
  domains over symmetric cones.\footnotetext{\emph{2000 Math Subject Classification:} 42B35, 32M15.}
\footnotetext{\emph{Keywords}: Bergman projection, symmetric cone.}}
\author{Cyrille Nana}
\address{Department of Mathematics, Faculty of Science,\\ University of Buea, P. O. Box 63,\\ Buea Cameroon}
\email{nana.cyrille@ubuea.cm}
\author{Beno\^it F. Sehba}
\address{Department of Mathematics, University of Ghana,\\ P. O. Box LG 62 Legon, Accra, Ghana}
\email{bfsehba@ug.edu.gh}

\begin{abstract}
We obtain some necessary and sufficient conditions for the boundedness of a family of positive operators defined on symmetric cones, we then deduce off-diagonal boundedness of associated Bergman-type operators in tube domains over symmetric cones.

\end{abstract}
 \maketitle

\section{Introduction}
\setcounter{equation}{0} \setcounter{footnote}{0}
\setcounter{figure}{0} 
The Bergman projection is the orthogonal projection from the Lebesgue space $L^2$ onto its closed subspace consisting of holomorphic functions. In this note, we obtain some answers on the question of the boundedness of the Bergman projection from the Lebesgue space $L^p$ to $L^q$ when $1\le p\le q\le \infty$ in the setting of tube domains over symmetric cones. For this, we first study a family of positive operators defined on symmetric cones. For our purpose, we introduce some definitions and notations.
\vskip .3cm 
Let $\O$ be an irreducible symmetric  cone 
in $\mathbb R^n$. 
There are a notion of rank and a notion 
of determinant function associated to each symmetric cone (see \cite{FK}). We write
 $r$ and  $\Delta$ for the rank and
the determinant function associated to $\O$ respectively. We recall that for $n\geq 3$, 
the Lorentz cone
$\Lambda_n$ of $\mathbb{R}^n$ is the set defined by
$$\Lambda_n=\{(y_1,\cdots,y_n)\in \mathbb R^n: y_1^2-\cdots-y_n^2>0,\,\,\,y_1>0\},$$
it is an example of symmetric cone of rank $r=2$ and its determinant function is given by the
Lorentz form
$$\Delta(y)=y_1^2-\cdots-y_n^2.$$  

\vskip .2cm
For $1\le p<\infty$ and $\nu \in \mathbb R$, we write
$L^p_\nu(\Omega)=L^p(\Omega,\Delta^{\nu-\frac{n}{r}}(y)dy)$
for the space of functions $f$ defined on  $\Omega$ that satisfy
 
$$\|f\|_{p,\nu}=||f||_{L^p_\nu(\Omega)}:=\left(\int_{\Omega}|f(y)|^p\Delta^{\nu-\frac{n}{r}}(y)dy\right)^{1/p}<\infty.$$
Let $\a,\beta,\ga$ be real parameters. We consider the integral operators $S=S_{\a,\beta,\ga}$
which are defined on $\O$ by
\Be\label{eq:conefamily}Sg(y)=\Da^\a(y)\int_\O\Da^{-\ga}(y+x)g(x)\Da^\beta(x)dx.\Ee
The above family can be seen as a generalization of the Hilbert-type operators considered in \cite{BanSeh}. It has been considered in several works related to the question of boundedness of Bergman operators (see for example \cite{BB2,BBGNPR,sehba} and the references therein).

In the first part of this paper, we consider the question of off-diagonal boundedness of the family $S_{\alpha,\beta,\gamma}$, that is its boundedness from  $L^p_\nu(\Omega)$ to $L^q_\mu(\Omega)$ for $1\le p\le q\le \infty$, $\nu,\mu\in \mathbb{R}$.

One application of the study of the family $S_{\alpha,\beta,\gamma}$ (see for example \cite{BBGNPR,sehba}) is that its boundedness implies the boundedness of a related family of Bergman-type operators defined below.
\vskip .2cm
For $\O$ an irreducible symmetric  cone 
in $\mathbb R^n$, we write $\mathcal D=\mathbb R^n+i\Omega$ for the tube domain over the  cone $\Omega$. Let us denote by $\mathcal{H}(\mathcal D)$ the space
of holomorphic functions on $\mathcal D$. For $1\le p<\infty$ and $\nu \in \mathbb R$, we write
$L^p_\nu(\mathcal D)=L^p(\mathcal
D,\Delta^{\nu-\frac{n}{r}}(y)dx\,dy)$
for the space of functions $f$ in $\mathcal D$
satisfying 
$$\|f\|_{p,\nu}=||f||_{L^p_\nu(\mathcal D)}:=\left(\int_{\mathcal D}|f(x+iy)|^p\Delta^{\nu-\frac{n}{r}}(y)dxdy\right)^{1/p}<\infty.$$
The weighted Bergman space is the set $A^p_\nu(\mathcal D)=L^p_\nu(\mathcal D)\cap \mathcal{H}(\mathcal D)$. It is a closed subspace of $L^p_\nu(\mathcal D)$ and it is well known that it is nontrivial only  for $\nu>\frac{n}{r}-1$ (see \cite{DD}). The usual Bergman space
$A^p(\mathcal D)$ corresponds to the case $\nu=\frac{n}{r}$.
\vskip .2cm
We also consider the mixed norm version of the above spaces. For $1\le p,q\leq \infty$, we write
$L^{p,q}_\nu(\mathcal D)=L^{p,q}(\mathcal
D,\Delta^{\nu-\frac{n}{r}}(y)dx\,dy)$
for the space of functions $f$ in $\mathcal D$
satisfying 
$$\|f\|_{p,q,\nu}=||f||_{L^{p,q}_\nu(\mathcal D)}:=\left(\int_{\O}\left(\int_{\mathbb{R}^n}|f(x+iy)|^pdx\right)^{q/p}\Delta^{\nu-\frac{n}{r}}(y)dy\right)^{1/q}<\infty$$
for $1\le q<\infty$, and 
$$\|f\|_{p,\infty}=||f||_{L^{p,\infty}(\mathcal D)}:=\sup_{y\in \O}\left(\int_{\mathbb{R}^n}|f(x+iy)|^pdx\right)^{1/p}<\infty.$$
The mixed norm Bergman space $A^{p,q}_\nu(\mathcal D)$ is defined as above. Note that the space $A^{p,\infty}(\mathcal D)$ is the Hardy space $\mathcal{H}^p(\mathcal{D})$.
\vskip .3cm
The Bergman projection $P_\nu$ is the orthogonal projection mapping $L^2_\nu(\mathcal
D)$ onto $A^2_\nu(\mathcal D)$. It is defined by
\Be P_{\nu}f(z)=\int_{\mathcal{D}}B_{\nu}(z,
w)f(w)dV_\nu(w) \Ee where
\Be\label{bergker} B_{\nu}(z,
w)=d_{\nu}\Da^{-\nu-\frac{n}{r}}(\frac{z-\overline w}{i}) \Ee is
the weighted Bergman kernel and $dV_\nu(w)=\Da^{\nu-\frac{n}{r}}(\Im w)dV(w)$, $dV$ being the Lebesgue measure on
$\mathbb{C}^n$(see \cite{BBGNPR}).
The question of whether
$P_\nu$ extends or not as a bounded operator on $L^p_\nu(\mathcal D)$ for $p\neq 2$ has attracted a
lot of attention in recent years (see \cite{BB1,BB2,BBGNPR,BBGR,BBGRS,B,BN,DD} and the references therein). 
Only the case of the tube domains over the Lorentz cone has been completely settled after the works \cite{BBGR,BBGRS,BBPR,GS} 
(cf. also \cite{BBGNPR} and \cite{B}), and the recent proof of the $l^2$-decoupling conjecture by Bourgain and 
Demeter \cite{BD} used as a key tool in \cite{BN} and \cite{BGN}; note that authors of \cite{BGN} present the complete answer to
this question for the case of the Pyateckii-Shapiro domain, which is
a Siegel domain of type II associated to the spherical cone of $\mathbb{R}^3.$
 
\vskip .3cm
In the second part of this work, we consider the same question but on the off-diagonal side for a family of
operators generalizing the Bergman projection.  
This family is
given by the integral operators $T=T_{\a,\beta,\ga}$ and $T^{+}=T_{\a,\beta,\ga}^+$ defined for functions in $C_c^{\infty}(\mathcal{D})$ 
by the
formula

$$Tf(z)=\Da^{\alpha}(\Im z)\int_{\mathcal{D}}B_{\ga}(z,w)f(w)\Da^{\ba}(\Im w)dV(w)$$
and
$$T^{+}f(z)=\Da^{\a}(\Im z)\int_{\mathcal{D}}|B_{\ga}(z,w)|f(w)\Da^{\beta}(\Im w)dV(w).$$
We note that the boundedness of $T^+$ implies the boundedness of $T$ and the reverse is not necessarily true.
It is clear that estimates obtained for this family of operators imply corresponding estimates for the weighted Bergman projection since 
$P_\nu=T_{0,\nu-\frac nr,\nu}$ for all $\nu>\frac nr-1.$
The $L^{p,q}_\nu(\mathcal D)$-boundedness of this family of operators has been considered in \cite{BBGR} for the case $T_{0,\mu-\frac nr,\mu}$ and
in \cite{BB2} for $T_{0,\mu-\frac nr,\mu+m}$ when $\O$ is the light cone.  Also, the second author in \cite[Theorem 1.1]{sehba} considered 
this family and obtained $L^{p,q}_\nu(\mathcal D)$ optimal results for the operator $T^+.$ In these works, was used the Schur Lemma with test functions 
which are generalized powers of the determinant function. Moreover, authors of \cite{BGN} extended some results of \cite{BBPR}
to the case of homogeneous Siegel domains of $\C^n$ of type II. They obtained sufficient conditions for the operator
$T^+_{\a,\mu-\tau,\mu+\a}$ in the vectorial weighted case, using the Schur Lemma. As was indicated in the introduction of \cite{Nana}
and explicitly presented in \cite{NT}, it is not yet known in general, with the techniques developed so far whether these sufficient 
conditions are necessary. However, with some restrictions on the weight or the choice of the cone (the Vinberg cone for example \cite{BNa}),
these sufficient conditions are necessary. 

We note that in \cite{sehba}, the estimates for the positive operators $T_{\alpha,\beta,\gamma}$ were applied in the characterization of the dual space of a Bergman space in the case where the associated Bergman projection is not necessarily bounded. Also, authors of \cite{BBGRS} made a heavy use of the family
$T_{\a,\nu-\frac nr,\nu+\a}$ in their characterization of Besov spaces of tube domains over symmetric cones.

 The aim of this paper is to extend to the setting of the tube domains over 
symmetric cones, the $L^p-L^q$ estimates of the Bergman projector, when $1\le p\le q\le \infty,$ considered by authors of
\cite{BanSeh}, who studied the case of upper half plane.  This leads us to more general result than those of \cite{sehba}. Our motivation is actually to know whether the Okikiolu test \cite{Oki},
which is a generalization of Schur test, used so far in various papers, is extendable in higher ranks. As we said earlier, it has been noticed recently that
in vectorial weighted cases, it is not clear that the Schur Lemma gives sharp results (see \cite{DD},\cite{NT} for example). 
The techniques that we develop here, do not impose the use of the generalized powers of the determinant function as
 test functions. Actually, our results lead to necessary and sufficient conditions for the $L^{p,q}(\mathcal D)$-boundedness of the Bergman operator 
 with positive Bergman kernel $P_\nu^+$ (see Corollary \ref{bounded} in the text), that coincide with classical ones obtained by 
 the Schur Lemma via generalized powers of the determinant function as
 test functions.
 
This work is divided in 5 sections. In Section 2, we state our results; in Section 3, we recall some of the key tools needed to establish our results;
in Section 4, we prove off-diagonal estimates for the operator $S_{\a,\ba,\ga}$ and the proofs of the off-diagonal estimates for the 
family of Bergman-type operators $T_{\a,\ba,\ga}$ are given in Section 5.

\vskip .2cm
\section{Statement of the results}
\subsection{Positive operators on the cone}
The following result provides the right relations between the parameters under which the operators $S_{\alpha,\beta,\gamma}$ are bounded.
\btheo\label{thm:main1cone}
Let $\nu,\mu\in\mathbb{R}$ and $1< p\leq q<\infty.$ Assume that $\frac \nu {p'}+\frac \mu q>0$. Then the following conditions are equivalent:

(a) The operator $S$ is bounded from $L^p_\nu(\O)$ into $L^q_\mu(\O).$

(b) The parameters satisfy
\Bea\label{0c}
\ga=\a+\beta+\frac nr-\frac \nu p+\frac \mu q,
\Eea
and 
\Bea\label{1c}
p\left(\beta-\ga+2\frac nr-1\right)-\frac nr+1<\nu<p(\beta+1)+\frac nr-1,
\Eea
and 
\Bea\label{2c}
\frac nr-1-q\a<\mu<q(\ga-\a)-\frac nr+1.
\Eea
\etheo
We also obtain the following first limit case.
\btheo\label{thm:main2cone}
Let $\nu,\mu\in\mathbb{R}$ and $1<  q<\infty.$ Assume that $\mu>0$. Then the following conditions are equivalent:
\begin{itemize}
\item[(a)] The operator $S$ is bounded from bounded from $L^1_\nu(\O)$ into $L^q_\mu(\O).$
\item[(b)] The parameters satisfy
\Bea\label{eq:cone21c}
\ga=\a+\beta+\frac nr-\nu +\frac \mu q
\Eea
and 
\Bea\label{eq:cone22c}
\ga>0
\Eea
and 
\Bea\label{eq:cone23c}
\frac nr-1-q\a<\mu<q(\ga-\a)-\frac nr+1.
\Eea
\end{itemize}

\etheo
The following second limit case is also proved.
\btheo\label{thm:main3cone}
Suppose $\nu\in \mathbb{R}$, and $1< p< \infty$. Then the
following conditions are equivalent:

\begin{itemize}
\item[(a)]

The operator $S$ is bounded from
$L_{\nu}^{p}(\O)$ to $L^{\infty}(\O).$

\item[(b)]

The parameters satisfy
\Bea
\ga=\a+\beta+\frac nr-\frac \nu p,
\Eea
and 
\Bea\label{eq:cone31c}
\nu<p(\beta+1)+\frac{n}{r}-1,
\Eea
and 
\Bea\label{eq:cone32c}
p(\a-\frac{n}{r}+1)>-\frac nr+1.
\Eea
\end{itemize}

\etheo
The necessity of the conditions in the above theorems is obtained using ideas from \cite{sehba}.  Our main tools for the sufficient parts will be a off-diagonal Schur's test due to G. O. Okikiolu and integrability conditions of the determinant function among others. We also refer to \cite{BanSeh,sehba1} for the corresponding one dimension results and applications.
\vskip .3cm
\subsection{Positive Bergman-type operators}
In the mixed norm case, following for example the proof of \cite[Theorem 1.1]{sehba} we obtain as consequence of Theorem \ref{thm:main1cone}, the following result.
\btheo\label{thm:main1}
Suppose $\nu\in \mathbb{R}$, $1< p< \infty$,and $1<q\le s<\infty$. Assume that $\frac \nu {p'}+\frac \mu q>0$. Then the
following conditions are equivalent:

\begin{itemize}
\item[(a)]

The operator $T_{\a,\beta,\ga}^+$ is bounded from
$L_{\nu}^{p,q}(\mathcal{D})$ to $L_{\mu}^{p,s}(\mathcal{D}).$

\item[(b)]

The parameters satisfy
\Bea
\ga=\a+\beta+\frac nr-\frac \nu q+\frac \mu s,
\Eea
and 
\Bea\label{eq:main11}
q\left(\beta-\ga+2\frac nr-1\right)-\frac nr+1<\nu<q(\beta+1)+\frac nr-1,
\Eea
and 
\Bea\label{eq:main12}
\frac nr-1-s\a<\mu<s(\ga-\a)-\frac nr+1.
\Eea
\end{itemize}

\etheo
In the same way, we have the following first limit case as consequence of Theorem \ref{thm:main2cone}.

\btheo\label{thm:main2}
Suppose $\nu\in \mathbb{R}$, and $1< p,s< \infty$. Assume that $\mu>0$. Then the
following conditions are equivalent:

\begin{itemize}
\item[(a)]

The operator $T_{\a,\beta,\ga}^+$ is bounded from
$L_{\nu}^{p,1}(\mathcal{D})$ to $L_{\mu}^{p,s}(\mathcal{D}).$

\item[(b)]

The parameters satisfy
\Bea
\ga=\a+\beta+\frac nr-\nu +\frac \mu s,
\Eea
and 
\Bea\label{eq:main21}
\ga>0,
\Eea
and 
\Bea\label{eq:main22}
\frac nr-1-s\a<\mu<s(\ga-\a)-\frac nr+1.
\Eea
\end{itemize}

\etheo
The following second limit case also follows from Theorem \ref{thm:main3cone}.
\btheo\label{thm:main3}
Suppose $\nu\in \mathbb{R}$, and $1< p,q< \infty$. Then the
following conditions are equivalent:

\begin{itemize}
\item[(a)]

The operator $T_{\a,\beta,\ga}^+$ is bounded from
$L_{\nu}^{p,q}(\mathcal{D})$ to $L^{p,\infty}(\mathcal{D}).$

\item[(b)]

The parameters satisfy
\Bea
\ga=\a+\beta+\frac nr-\frac \nu q,
\Eea
and 
\Bea\label{eq:main31}
\nu<q(\beta+1)+\frac{n}{r}-1,
\Eea
and 
\Bea\label{eq:main32}
q(\a-\frac{n}{r}+1)>-\frac nr+1.
\Eea
\end{itemize}

\etheo
\vskip .3cm
Finally, we have the following which gives sufficient conditions for the boundedness of operators $T_{\a,\beta,\ga}^+$ from $L_{\nu}^{p}(\mathcal{D})$ to $L_{\mu}^{q}(\mathcal{D})$.
\btheo\label{thm:main4}
Suppose $\nu\in \mathbb{R}$ and $1< p\le q< \infty$. Assume that $\frac{\nu+\frac nr}{p'}+\frac{\mu+\frac nr}{q}>0$, and that the parameters satisfy
\Be\label{eq:main40}
\ga=\a+\beta+\frac nr-\frac{1}{p}(\nu+\frac nr)+\frac{1}{q}(\mu+\frac nr),
\Ee

and 
\Bea\label{eq:main41}
\nu<p(\beta+1)+\left(\frac nr-1\right)\left(1-\frac{p}{q}\right),
\Eea
and 
\Bea\label{eq:main42}
-q\a+ \left(\frac nr-1\right)\left(1+\frac{q}{p'}\right)<\mu.
\Eea
Then the operator $T_{\a,\beta,\ga}^+$ is bounded from
$L_{\nu}^{p}(\mathcal{D})$ to $L_{\mu}^{q}(\mathcal{D})$.
\etheo

For the proof of Theorem \ref{thm:main4}, we directly appeal to the same off-diagonal Schur's test of G. O. Okikiolu.
\vskip .2cm
Once again, we recall that the boundedness of $T_{\a,\beta,\ga}^+$ implies the boundedness of $T_{\a,\beta,\ga}$ although the boundedness of $T_{\a,\beta,\ga}$ is expected in a larger range than the one of $T_{\a,\beta,\ga}^+$.
\vskip .2cm
Given two positive quantities $A$ and $B$, the notation $A\lesssim B$ (resp. $B\lesssim A$) will mean that there is an universal constant $C>0$ such that $A\le CB$ (resp. $B\le CA$). When $A\lesssim B$ and $B\lesssim A$, we write $A\backsimeq B$. Notation $C_\alpha$ means that the constant $C$ depends on the parameter $\alpha$.

\section{Preliminaries}
We give here some useful tools needed in our presentation. We refer to \cite{DD} for the following results. For all 
vector ${\bf s}=(s_1,\cdots,s_r),$ we denote 
$$\Da^{\bf s}(x)=\Da^{s_1-s_2}_1(x)\Da^{s_2-s_3}_2(x)\cdots\Da^{s_r}_r(x)$$
the generalized power of the determinant function. Here $\Da_j(x)$'s are the {\it principal minors} of $x\in\O$ with respect to the Jordan algebra
structure of the symmetric cone $\O.$ Also $\Da_r(x):=\Da(x),$ the determinant of $x.$

Put $(r-1)\frac{d}{2}=\frac{n}{r}-1$. We recall the following result.

\blem\label{lem:integralcone}
Let $v\in \O$ and ${\bf s,t} \in \mathbb C^r.$ The integral
$$\int_{\Omega} \Delta^{\bf s} (y+v)\Delta^{{\bf t}-\frac nr} (y)dy$$
converges if for every  $j=1,\cdots,r$ we have 
$\Re{t_j} >(r-j)\frac d2$ and $\Re ({ s_j +t_j}) <-(j-1)\frac d2.$ In this case, we have
$$\int_{\Omega} \Delta^{\bf s} (y+v)\Delta^{{\bf t}-\frac nr} (y)dy=C_{{\bf s}, {\bf t}} \Delta^{{\bf s}+{\bf t}} (v).$$
\elem
We also need the following result.
\blem\label{lem:Apqfunction} Let $\a$ be real. Then
the function
$f(z)=\Da^{-\a}(\frac{z+it}{i})$, with $t \in \O$, belongs to
$L_{\nu}^{p,q}(\mathcal{D})$ if and only if $\nu>\frac{n}{r}-1$ and $\a >
max\left(\frac{2\frac{n}{r} -1}{p},\frac{n}{rp} + \frac{\nu
+\frac{n}{r}-1}{q}\right)$. In this
case,$$||f||_{L_{\nu}^{p,q}}^q=C_{\a,p,q}\Da^{-q\a +\frac{nq}{rp}
+ \nu}(t).$$
\elem
\begin{proof} See \cite[Lemma 3.20]{BBGNPR}.
\end{proof}
The following extension of the Schur's test is due to G. O. Okikiolu \cite{Oki}.
\blem\label{lem:okikiolu}
Let $p,r,q$ be positive numbers such that $1<p\le r$ and $\frac{1}{p} + \frac{1}{q} = 1.$ Let $K(x,y)$ be a non-negative measurable function on $X\times Y$ and suppose there exist $0<t\le 1$, measurable functions $\phi_1:X\rightarrow (0,\infty),\quad\phi_2:Y\rightarrow (0,\infty)$ and nonnegative constants $M_1,M_2$ such that
\begin{eqnarray}
\label{ee1} 
\int_X K(x,y)^{tq}\phi_1^{q}(y)\mathrm{d}\mu(y) &\le & M_1^{q}\phi_2^{q}(x)\qquad\mbox{a.e on}\quad Y\quad\mbox{and}\\ 
\label{ee2}
\int_Y K(x,y)^{(1-t)r}\phi_2^r(x)\mathrm{d}\nu(x) &\le  &M_2^r\phi_1^r(y)\qquad\mbox{a.e on}\quad X.
\end{eqnarray} 
If $T$ is given by $$Tf(x)=\int_X\!\!\!f(y)K(x,y)\mathrm{d}\mu(y)$$ where $f\in L^p(X,\mathrm{d}\mu),$ then $T:L^p(X,\mathrm{d}\mu)\longrightarrow L^r(Y,\mathrm{d}\nu)$ is bounded and for each $f\in L^p(X,\mathrm{d}\mu)$, $$\left\|Tf\right\|_{L^r(Y,\mathrm{d}\nu)}\le M_1M_2\|f\|_{L^p(X,\mathrm{d}\mu)}.$$
\elem
The following limit case of Okikiolu result is proved in \cite{Zhao}.
\blem\label{lem:okikiolulimitcase}
Let $\mu$ and $\nu$ be positive measures on the space $X$ and let $K(x, y)$ be non-negative measurable functions on $X\times Y.$ Let $T$ be the integral operator with kernel $K(x,y)$ defined by $$Tf(x) = \int_X\!\! f(y)K(x, y)\mathrm{d}\mu(y).$$ Suppose $1=p\le q<\infty.$ Let $\gamma$ and $\delta$ be two real numbers such that $\gamma + \delta = 1.$ If there exist positive functions $h_1$ and $h_2$ with positive constants $C_1$ and $C_2$ such that $$\mbox{ess}\sup_{y\in Y}h_1(y)K(x,y)^{\gamma}\le C_1h_2(x)\quad\mbox{for almost all } x\in X\quad\mbox{and}$$ $$\int_X\!\!h_2(x)^qK(x,y)^{\delta q}\mathrm{d}\nu\le C_2h_1(y)^q\quad\mbox{for almost all } y\in Y,$$ then $T$ is bounded from $L^1(X,\mathrm{d}\nu)$ into $L^q(X,\mathrm{d}\nu)$ and the norm of this operator does not exceed $C_1C_2^{\frac{1}{q}}.$
\elem

\section{Off-diagonal estimates for the family $S_{\alpha,\beta,\gamma}$}
Let $\Omega$ be an irreducible symmetric cone in $\mathbb{R}^n$ of rank $r.$ Let $\a,\beta,\ga,\nu,\mu$ be real parameters. In this section, we obtain necessary and sufficient conditions for the boundedness of the operators $S=S_{\a,\beta,\ga}$ from $L_{\nu}^{p}(\O)$ to $L_{\mu}^{q}(\O)$, $1\le p\le q\le \infty$. We recall that the integral operators $S=S_{\a,\beta,\ga}$
are defined on $\O$ by
$$Sg(y)=\Da^\a(y)\int_\O\Da^{-\ga}(y+x)g(x)\Da^\beta(x)dx.$$
This operator has already been considered by the author of \cite{sehba}. He found necessary and sufficient conditions under which this operator is bounded on $L_{\nu}^{q}(\O)$, using the Schur Lemma. He has applied this lemma with vectorial powers of the determinant function. We shall prove here, a more general result, an 
off-diagonal estimate of this operator. It will require the use of a more general Schur's test result for sufficiency: the Okikiolu's lemma. 

\subsection{Sufficiency for the boundedness of $S_{\alpha,\beta,\gamma}$. }
Let us start by proving the following result.
\blem\label{lem:conesuff1}
Let $\nu,\mu\in\mathbb{R}$ and $1< p\leq q<\infty.$ Assume that $\frac \nu {p'}+\frac \mu q>0$. Suppose that the parameters satisfy
\Bea\label{0}
\ga=\a+\beta+\frac nr-\frac \nu p+\frac \mu q,
\Eea
and 
\Bea\label{1}
p\left(\beta-\ga+2\frac nr-1\right)-\frac nr+1<\nu<p(\beta+1)+\frac nr-1,
\Eea
and 
\Bea\label{2}
\frac nr-1-q\a<\mu<q(\ga-\a)-\frac nr+1.
\Eea
Then the operator $S$ is bounded from $L^p_\nu(\O)$ into $L^q_\mu(\O).$
\elem



\begin{proof}
Assume $$\ga=\a+\beta+\frac nr-\frac \nu p+\frac \mu q$$ and let
$$\omega=\a+\beta-\ga-\nu+\frac nr=-\left(\frac \nu {p'}+\frac \mu q\right)<0.$$
From the left inequality of (\ref{1}), we have $p(\beta-\ga+2\frac nr-1)-\frac nr+1<\nu$ which is equivalent 
to $\frac{\nu+\frac nr-1}{p}>\beta-\ga+2\frac nr-1.$ Multiplying this last 
inequality by $\omega<0$ yields $\frac{\nu+\frac nr-1}{p}\omega<(\beta-\ga+2\frac nr-1)\omega$ i.e.
\Beas
\frac{\nu+\frac nr-1}{p}\omega+\frac{\beta-\ga+2\frac nr-1}{p'}\nu<-\frac{\beta-\ga+2\frac nr-1}{q}\mu.
\Eeas
This is equivalent to
\Bea\label{1.3}
\frac{\nu+\frac nr-1}{p'}\omega-\frac{\beta-\ga-\nu+\frac nr}{p'}\nu>\frac{\beta-\ga-\nu+\frac nr}{q}\mu.
\Eea

From the right inequality of (\ref{1}), we have $\nu<p(\beta+1)+\frac nr-1$ which is equivalent to $\frac{\nu-\frac nr+1}{p}<\beta+1.$ Hence,
$\nu-\frac nr+1-\frac 1{p'}(\nu-\frac nr+1)<\beta+1$ so that $\beta-\nu+\frac nr+\frac 1{p'}(\nu-\frac nr+1)>0.$ Multiplying this last 
inequality by $\omega$ yields $(\beta-\nu+\frac nr)\omega+\frac 1{p'}(\nu-\frac nr+1)\omega<0$, i.e.
\Bea\label{1.1}
\frac{\nu-\frac nr+1}{p'}\omega-\frac{\beta-\nu+\frac nr}{p'}\nu<\frac{\beta-\nu+\frac nr}{q}\mu.
\Eea

From the left inequality of (\ref{2}), we have $\frac nr-1-q\a<\mu$ which is equivalent 
to $-\frac{\mu-\frac nr+1}{q}-\a<0.$ Multiplying this last 
inequality by $\omega$ yields $-\frac{\mu-\frac nr+1}{q}\omega-\a w>0$, i.e.
\Bea\label{1.2}
\frac{\mu-\frac nr+1}{q}\omega-\a\frac{\mu}{q}<\a\frac{\nu}{p'}.
\Eea

From the right inequality of (\ref{2}), we have $\mu<q(\ga-\a)-\frac nr+1$ which is equivalent to $\frac{\mu+\frac nr-1}{q}<\ga-\a.$ 
Multiplying this last 
inequality by $\omega$ yields $\frac{\mu+\frac nr-1}{q}\omega>(\ga-\a)\omega$, i.e.
\Bea\label{1.4}
\frac{\mu+\frac nr-1}{q}\omega+\frac{\ga-\a}{q}\mu>-\frac{\ga-\a}{p'}\nu.
\Eea

The inequalities (\ref{1.1}), (\ref{1.2}), (\ref{1.3}) and (\ref{1.4}) yield the existence of two real numbers $u$ and $v$ such that
\Bea\label{u}
\left\{\Ba{l}\frac{\nu-\frac nr+1}{p'}\omega-\frac{\beta-\nu+\frac nr}{p'}\nu<u\omega+(\beta-\nu+\frac nr)(v-u)<\frac{\beta-\nu+\frac nr}{q}\mu\\\\
\frac{\mu-\frac nr+1}{q}\omega-\a\frac{\mu}{q}<v\omega+\a(u-v)<\a\frac{\nu}{p'}\Ea\right.
\Eea
and
\Bea\label{v}
\left\{\Ba{l}\frac{\nu+\frac nr-1}{p'}\omega-\frac{\beta-\ga-\nu+\frac nr}{p'}\nu>u\omega+(\beta-\ga-\nu+\frac nr)(v-u)>\frac{\beta-\ga-\nu+\frac nr}{q}\mu\\\\
\frac{\mu+\frac nr-1}{q}\omega+\frac{\ga-\a}{q}\mu>v\omega+(\ga-\a)(v-u)>-\frac{\ga-\a}{p'}\nu.\Ea\right.
\Eea
Now, (\ref{u}) is equivalent to
\Bea\label{u1}
\left\{\Ba{l}-\frac{\beta-\nu+\frac nr}{\omega}\left[-\frac{\mu}{q}-u+v\right]<u<\frac{\nu-\frac nr+1}{p'}+\frac{\beta-\nu+\frac nr}{\omega}\left[-\frac{\nu}{p'}+u-v\right]\\\\
-\frac{\a}{\omega}\left[-\frac{\nu}{p'}+u-v\right]<v<\frac{\mu-\frac nr+1}{q}+\frac{\a}{\omega}\left[-\frac{\mu}{q}-u+v\right]\Ea\right.
\Eea
and (\ref{v}) is equivalent to
\Bea\label{v1}
\left\{\Ba{l}\frac{\nu+\frac nr-1}{p'}+\frac{\beta-\ga-\nu+\frac nr}{\omega}\left[-\frac{\nu}{p'}+u-v\right]<u<-\frac{\beta-\ga-\nu+\frac nr}{\omega}\left[-\frac{\mu}{q}-u+v\right]\\\\
\frac{\mu+\frac nr-1}{q}-\frac{\ga-\a}{\omega}\left[-\frac{\mu}{q}-u+v\right]<v<\frac{\ga-\a}{\omega}\left[-\frac{\nu}{p'}+u-v\right].\Ea\right.
\Eea
Let
$$t=\frac{-\frac{\nu}{p'}+u-v}{\omega}$$
then
$$1-t=\frac{-\frac{\mu}{q}-u+v}{w}.$$
Since $\omega<0,$ we choose $u$ and $v$ such that $0<v-u<\frac{\mu}{q}.$ Thus, we have $0<t<1.$
Therefore (\ref{u1}) and (\ref{v1}) become
\Bea\label{u2}
\left\{\Ba{l}-(\beta-\nu+\frac nr)(1-t)<u<\frac{\nu-\frac nr+1}{p'}+(\beta-\nu+\frac nr)t\\\\
-\a t<v<\frac{\mu-\frac nr+1}{q}+\a(1-t)\Ea\right.
\Eea
and 
\Bea\label{v2}
\left\{\Ba{l}\frac{\nu+\frac nr-1}{p'}+(\beta-\ga-\nu+\frac nr)t<u<-(\beta-\ga-\nu+\frac nr)(1-t)\\\\
\frac{\mu+\frac nr-1}{q}-(\ga-\a)(1-t)<v<(\ga-\a)t\Ea\right.
\Eea
respectively.

We shall now use the Okikiolu test to conclude. To this effect, we observe that the kernel of the operator $S:L^p_\nu(\O)\to L^q_\mu(\O)$ with respect to the measure $\Da^{\nu-\frac nr}(x)dx$ is given by
$$K(y,x)=\Da^\a(y)\Da^{-\ga}(y+x)\Da^{\beta-\nu+\frac nr}(x).$$
Consider the positive functions $\phi_1(x)=\Da^{-u}(x)$ and $\phi_2(y)=\Da^{-v}(y).$ Then
\Beas
I_1&=&\int_\O K(y,x)^{tp'}\phi_1(x)^{p'}\Da^{\nu-\frac nr}(x)dx\\
&=&\Da^{tp'\a}(y)\int_\O\Da^{-tp'\ga}(y+x)\Da^{tp'(\beta-\nu+\frac nr)-p'u+\nu-\frac nr}(x)dx.
\Eeas
The last integral above converges because from the right inequality in (\ref{u2}) involving $u$ and 
the left inequality in (\ref{v2}) involving $u$, we have $tp'(\beta-\nu+\frac nr)-p'u+\nu>\frac nr-1$ and 
$-tp'\ga+tp'(\beta-\nu+\frac nr)-p'u+\nu<-\frac nr+1$ respectively. It follows using Lemma \ref{lem:integralcone} that
\Beas
I_1&=&\Da^{tp'\a}(y)\int_\O\Da^{-tp'\ga}(y+x)\Da^{tp'(\beta-\nu+\frac nr)-p'u+\nu-\frac nr}(x)dx\\
&=&C_1\Da^{tp'\a-tp'\ga+tp'(\beta-\nu+\frac nr)-p'u+\nu}(y)\\
&=&C_1\Da^{-p'v}(y)=C_1\phi_2(y)^{p'}
\Eeas
since $tp'\a-tp'\ga+tp'(\beta-\nu+\frac nr)-p'u+\nu=tp'w-p'u+\nu=p'(-\frac{\nu}{p'}+u-v)-p'u+\nu=-p'v.$

Now,
\Beas
I_2&=&\int_\O K(y,x)^{(1-t)q}\phi_2(y)^{q}d\mu(y)\\
&=&\Da^{(1-t)q(\beta-\nu+\frac nr)}(x)\int_\O\Da^{-(1-t)q\ga}(y+x)\Da^{(1-t)q\a-qv+\mu-\frac nr}(y)dy.
\Eeas
The last integral above converges because from the right inequality in (\ref{u2}) involving $v$ and 
the left inequality in (\ref{v2}) involving $v$, we have $$(1-t)q\a-qv+\mu>\frac nr-1\,\,\, \textrm{and}\\ 
-(1-t)q\ga+(1-t)q\a-qv+\mu<-\frac nr+1$$ respectively. It follows using Lemma \ref{lem:integralcone} again that
\Beas
I_2&=&\Da^{(1-t)q(\beta-\nu+\frac nr)}(x)\int_\O\Da^{-(1-t)q\ga}(y+x)\Da^{(1-t)q\a-qv+\mu-\frac nr}(y)dy\\
&=&C_2\Da^{(1-t)q(\beta-\nu+\frac nr)-(1-t)q\ga+(1-t)q\a-qv+\mu}(x)\\
&=&C_2\Da^{-qu}(x)=C_2\phi_1(x)^{q}
\Eeas
since 
\Beas
&& (1-t)q(\beta-\nu+\frac nr)-(1-t)q\ga+(1-t)q\a-qv+\mu\\ &=& (1-t)q\omega-qv+\mu\\ &=& q\omega-qt\omega-qv+\mu=-q\frac{\nu}{p'}-\mu-q(\frac{-\nu}{p'}+u-v)-qv+\mu\\ &=&
-qu.
\Eeas 
Thus by the Okikiolu test, we conclude that $S:L^p_\nu(\O)\to L^q_\mu(\O)$ is bounded.
\end{proof}
\begin{remark}
One could have also taken the powers $u$ and $v$ in our test functions as vectors, i.e. ${\bf u}=(u_1,\ldots,u_r)$ and ${\bf v}=(v_1,\ldots,v_r)$. Indeed, in doing this, Lemma \ref{lem:integralcone} again shows that for every $j=1,\ldots,r$, we should have
\Beas
\left\{\Ba{l}tp'(\ba-\nu+\frac nr)-p'u_j+\nu>(r-j)\frac{d}{2}\\-tp'\ga+tp'(\ba-\nu+\frac nr)-p'u_j+\nu<-(j-1)\frac d2\Ea\right.
\Eeas
and 

\Beas
\left\{\Ba{l}(1-t)q\a-qv_j+\mu>(r-j)\frac{d}{2}\\-(1-t)q\ga+(1-t)q\a-qv_j+\mu<-(j-1)\frac d2\Ea.\right.
\Eeas
But an observation of the systems (\ref{u2}) and (\ref{v2}) shows that $\gamma>0$ so that the interval $(\frac{\nu+\frac nr-1}{p'}+(\beta-\ga-\nu+\frac nr)t,\frac{\nu-\frac nr+1}{p'}+(\beta-\nu+\frac nr)t)$ is not empty; we also see that the interval $(\frac{\mu+\frac nr-1}{q}-(\ga-\a)(1-t),\frac{\mu-\frac nr+1}{q}+\a(1-t))$ is not empty. These facts clearly secure the existence of the vectors ${\bf u}$ and ${\bf v}$. 
\end{remark}
The following provides sufficient conditions for the boundedness of $S_{\alpha,\beta,\gamma}$ from $L^1_\nu(\O)$ into $L^q_\mu(\O)$.
\blem\label{lem:main2cone}
Let $\nu,\mu\in\mathbb{R}$ and $1<  q<\infty.$ Assume that $\mu>0$. Suppose that
the parameters satisfy
\Bea\label{eq:cone21}
\ga=\a+\beta+\frac nr-\nu +\frac \mu q
\Eea
and 
\Bea\label{eq:cone22}
\ga>0
\Eea
and 
\Bea\label{eq:cone23}
\frac nr-1-q\a<\mu<q(\ga-\a)-\frac nr+1.
\Eea
 Then the operator $S$ is bounded from $L^1_\nu(\O)$ into $L^q_\mu(\O).$

\elem
\begin{proof}
Assume that $\ga=\a+\beta+\frac nr-\nu +\frac \mu q$ and the parameters satisfy (\ref{eq:cone22}) and (\ref{eq:cone23}). Note that combining (\ref{eq:cone21}) and (\ref{eq:cone23}), one obtains that 
$$\beta-\gamma+\frac nr<\nu<\beta+\frac nr$$
which corresponds to (\ref{1}) for $p=1$.
Using the notations of the proof of the previous result, we now have $$\omega=\a+\beta+\frac nr-\nu-\ga=-\frac \mu q<0,$$
and $$t=\frac{u-v}{\omega}.$$ Now (\ref{u2}) reduces to
\Bea\label{u22}
\left\{\Ba{l}-(\beta-\nu+\frac nr)(1-t)<u<(\beta-\nu+\frac nr)t\\\\
-\a t<v<\frac{\mu-\frac nr+1}{q}+\a(1-t)\Ea\right.
\Eea
and (\ref{v2}) reduces to
\Bea\label{v22}
\left\{\Ba{l}(\beta-\ga-\nu+\frac nr)t<u<-(\beta-\ga-\nu+\frac nr)(1-t)\\\\
\frac{\mu+\frac nr-1}{q}-(\ga-\a)(1-t)<v<(\ga-\a)t.\Ea\right.
\Eea
The second condition in Lemma \ref{lem:okikiolulimitcase} can be then checked as in the proof of the previous theorem. 
Let us also check that the other test condition in the same lemma is verified. That is there are a constant $C>0$ and positive 
functions $\phi_1$ and $\phi_2$ such that 
$$\sup_{x\in \O}\,\phi_1(x)K(y,x)^t\le C\phi_2(y)\,\,\,\textrm{for almost every}\,\,\,y\in \O.$$
Our positive functions are still given by $\phi_1(x)=\Da^{-u}(x)$ and $\phi_2(y)=\Da^{-v}(y)$. 
We then have for every $x\in\O,$
$$\phi_1(x)K(y,x)^t\phi_2(y)^{-1}=\Da^{-u+t\left(\beta-\nu+\frac nr\right)}(x)\Da^{-\ga t}(y+x)\Da^{v+t\alpha}(y).$$
Observe that from (\ref{u22}), $v+t\alpha>0.$ Hence using that for all $a,b\in \O$, $$\Da(a+b)\geq \Da(a)$$ (see \cite{FK}), we obtain that
$\Da^{v+t\alpha}(y)\leq \Da^{v+t\alpha}(y+x)$ so that
\Bea
\phi_1(x)K(y,x)^t\phi_2(y)^{-1}\leq\Da^{-u+t\left(\beta-\nu+\frac nr\right)}(x)\Da^{-\ga t+v+t\alpha}(y+x).
\Eea
Moreover, from (\ref{u22}), we have $-u+t\left(\beta-\nu+\frac nr\right)>0.$ Therefore, 
$$\Da^{-u+t\left(\beta-\nu+\frac nr\right)}(x)\leq\Da^{-u+t\left(\beta-\nu+\frac nr\right)}(y+x).$$ It follows that for every $x\in\O,$
\Bea
\phi_1(x)K(y,x)^t\phi_2(y)^{-1}\leq\Da^{-u+t\left(\beta-\nu+\frac nr\right)-\ga t+v+t\alpha}(y+x)=1
\Eea
for almost every $y\in\O,$ since
$$-u+t\left(\beta-\nu+\frac nr\right)-\ga t+v+t\alpha=-u+v+t\omega=-u+v+u-v=0.$$
\end{proof}
We now prove the following.
\blem\label{lem:main3cone}
Let $\nu\in \mathbb{R}$, and $1< p< \infty$. Suppose that the parameters satisfy
\Bea
\ga=\a+\beta+\frac nr-\frac \nu p,
\Eea
and 
\Bea\label{eq:cone31}
\nu<p(\beta+1)+\frac{n}{r}-1,
\Eea
and 
\Bea\label{eq:cone32}
p(\a-\frac{n}{r}+1)>-\frac nr+1.
\Eea
Then the operator $S$ is bounded from $L_\nu^p(\O)$ to $L^\infty(\O)$.

\elem
\begin{proof}
 Let $g\in L_\nu^p(\O)$. Assume that $\ga=\a+\beta+\frac nr-\frac \nu p$. Then using the H\"older's inequality and Lemma \ref{lem:integralcone}, we obtain
\Beas
|Sg(y)| &=& \left|\Da^{\a}(y)\int_\O\Da^{-\ga}(y+x)g(x)\Da^\beta(x)dx\right|\\
&=& \left|\Da^{\a}(y)\int_\O\Da^{-\ga}(y+x)g(x)\Da^{\beta-\nu+\frac nr}(x)dV_\nu(x)\right|\\ 
&\le& \Da^{\a}(y)\int_\O\Da^{-\ga}(y+x)|g(x)|\Da^{\beta-\nu+\frac nr}(x)dV_\nu(x)\\
&\le& \|g\|_{p,\nu}\Da^{\a}(y)\left(\int_\O\Da^{-p'\ga}(y+x)\Da^{p'(\beta-\nu+\frac nr)+\nu-\frac nr}(x)dx\right)^{1/p'}\\ 
&=& C\|g\|_{p,\nu}\Da^{\a-\ga+\beta-\nu+\frac nr+\frac{\nu}{p'}}(y)\\ &=& C\|g\|_{p,\nu}.
\Eeas 
Hence $$\sup_{y\in \O}|Sg(y)|\le C\|g\|_{p,\nu}.$$ The proof is complete.
\end{proof}
\subsection{Necessity for the boundedness of $S_{\alpha,\beta,\gamma}$. }
\blem\label{lem:main1coneness}
Let $\nu,\mu\in\mathbb{R}$ and $1< p\leq q<\infty.$ Assume that $\frac \nu {p'}+\frac \mu q>0$. If the operator $S$ is bounded from $L^p_\nu(\O)$ into $L^q_\mu(\O)$, then the parameters satisfy the conditions (\ref{0}), (\ref{1}) and (\ref{2}).
\elem
\begin{proof}
 We start with the proof of the homogeneity condition (\ref{0}), that is, 
$$\ga=\a+\beta+\frac nr-\frac \nu p+\frac \mu q.$$
For this, we recall that the determinant function is homogeneous of degree $r$ (see \cite{FK}). Let $R>0$. 
To any $f\in L^p(\O, \Da^{\nu-\frac nr}(y)dy)$, we associate the function $f_R$ defined by $f_R(y)=f(Ry)$. One easily checks that 
$$\|f_R\|_{p,\nu}=R^{-r\frac{\nu}{p}}\|f\|_{p,\nu}.$$
An easy change of variable combined with the fact that the determinant function is homogeneous of degree $r$ provides
$$S(f_R)(y)=R^{-n+r(\gamma-\beta-\alpha)}Sf(Ry).$$
It follows using again the homogeneity of the determinant function that 
$$\|Sf_R\|_{q,\mu}=R^{-n+r(\gamma-\beta-\alpha-\frac{\mu}{q})}\|Sf\|_{q,\mu}.$$
From the boundedness of the operator $S$, we have that there is a constant $C$ such that for any $f\in L^p(\O, \Da^{\nu-\frac nr}(y)dy)$,
\Beas R^{-n+r(\gamma-\beta-\alpha-\frac{\mu}{q})}\|Sf\|_{q,\mu}=\|Sf_R\|_{q,\mu} &\le& C\|f_R\|_{1,\nu} =
CR^{-r\frac{\nu}{p}}\|f\|_{p,\nu},\Eeas
that is, 
$$ R^{-n+r(\gamma-\beta-\alpha-\frac{\mu}{q}+\frac{\nu}{p})}\|Sf\|_{q,\mu}\le C\|f\|_{p,\nu}.$$
As the latter holds for every $f\in L^p(\O, \Da^{\nu-\frac nr}(y)dy)$ and as $R$ was taken arbitrary, we should necessarily have
$$-n+r(\gamma-\beta-\alpha-\frac{\mu}{q}+\frac{\nu}{p})=0,$$
which leads to $\gamma=\alpha+\beta+\frac nr-\frac{\nu}{p}+\frac{\mu}{q}$. This proves (\ref{0}).
\vskip .2cm
It is well known that the symmetric cone $\O$ induces in $V\equiv \mathbb{R}^n$ a structure of Euclidean
Jordan algebra. We denote by $\underline{\mathbf e}$ the identity element in $V$. Let $g=\chi_{B(\underline{\mathbf e},1)}$, where $B(\underline{\mathbf e},1)$ is the Euclidean ball about $\underline{\mathbf e}$ with radius $1$. Following \cite{sehba}, we have $$Sg(y)\simeq \Da^\alpha(y)\Da^{-\gamma}(y+\underline{\mathbf e}).$$
It follows that if $S$ is bounded from $L^p(\O, \Da^{\nu-\frac nr}(y)dy)$ to $L^q(\O, \Da^{\mu-\frac nr}(y)dy)$, 
then the function $\Da^\alpha(y)\Da^{-\gamma}(y+\underline{\mathbf e})$ is in $L^q(\O, \Da^{\mu-\frac nr}(y)dy)$, 
which means that $$\int_{\O}\Da^{q\alpha+\mu-\frac nr}(y)\Da^{-q\gamma}(y+{\bf e})dy<\infty.$$
It follows from Lemma \ref{lem:integralcone} that we should have $$q\alpha+\mu>\frac nr-1\quad\mbox{and}\quad -q\gamma+q\alpha+\mu<-\frac nr+1.$$ That is $$\frac nr-1-q\a<\mu<q(\ga-\a)-\frac nr+1$$ which is condition (\ref{2}).
To prove the necessity of the condition (\ref{1}),
we proceed by duality. We have that the boundedness of $S$ from $L^p(\O, \Da^{\nu-\frac nr}(y)dy)$ to $L^q(\O, \Da^{\mu-\frac nr}(y)dy)$
is equivalent to the boundedness of the adjoint $S^*$ of $S$ from $L^{q'}(\O, \Da^{\mu+q'(\nu-\mu)-\frac nr}(y)dy)$ to $L^{p'}(\O, \Da^{\nu-\frac nr}(y)dy)$. We note that $S^*$ is given by 
$$S^*g(v)=\Da^{\beta-\nu+\frac nr}(v)\int_{\O}\Da^{-\gamma}(y+v)g(y)\Da^{\alpha+\nu-\frac nr}(y)dy.$$
Proceeding as above, we obtain that the function $\Da^{\beta-\nu+\frac nr}(v)\Da^{-\gamma}(v+\underline{\mathbf e})$
must belong to $L^{p'}(\O, \Da^{\nu-\frac nr}(y)dy)$. Using again Lemma \ref{lem:integralcone}, we see that we must have
$(\ba - \nu + \frac{n}{r})p'+ \nu>\frac{n}{r} -1$ and $-p'\ga +
(\ba - \nu + \frac{n}{r})p'+ \nu < -\frac{n}{r} +
1$, which is equivalent to $\nu < p(\ba + 1) + \frac{n}{r} -1$ and
$\nu > p(\beta -\ga  + 2\frac{n}{r} -1) -\frac{n}{r} + 1$. This
completes the proof of the lemma.

\end{proof}
Let us also prove the necessity of the conditions in Lemma \ref{lem:main2cone}.
\blem\label{lem:main2coneness}
Let $\nu,\mu\in\mathbb{R}$ and $1<  q<\infty.$ Assume that $\mu>0$. If the operator $S$ is bounded from $L^1_\nu(\O)$ into $L^q_\mu(\O)$, then the parameters satisfy the conditions (\ref{eq:cone21}), (\ref{eq:cone22}) and (\ref{eq:cone23}).
\elem
\begin{proof}
The necessity of the homogeneity condition
$$\gamma=\alpha+\beta+\frac nr-\nu+\frac{\mu}{q}$$ follows as in the proof of (\ref{0}).
The necessity of the condition (\ref{eq:cone23}) follows by taking $p=1$ in the proof of the inequality (\ref{2}) above. 
To prove the necessity of the condition (\ref{eq:cone22}), we proceed again by duality. We have that the boundedness of $S$ 
from $L^1(\O, \Da^{\nu-\frac nr}(y)dy)$ to $L^q(\O, \Da^{\mu-\frac nr}(y)dy)$ is equivalent to the boundedness of the adjoint $S^*$ of $S$
from $L^{q'}(\O, \Da^{\mu+q'(\nu-\mu)-\frac nr}(y)dy)$ to $L^\infty(\O)$. We have seen that $S^*$ is given by 
$$S^*g(v)=\Da^{\beta-\nu+\frac nr}(v)\int_{\O}\Da^{-\gamma}(y+v)g(y)\Da^{\alpha+\nu-\frac nr}(y)dy.$$
Testing again with the function $g=\chi_{B(\underline{\mathbf e},1)}$, where $B(\underline{\mathbf e},1)$ is the Euclidean ball about $\underline{\mathbf e}$ with radius $1$, we obtain that the function $\Da^{\beta-\nu+\frac nr}(v)\Da^{-\gamma}(v+\underline{\mathbf e})$
should belong to $L^\infty(\O)$. That is for any $v\in \O$, 
\begin{equation}\label{eq:nessnotdirect}
\Da^{\beta-\nu+\frac nr}(v)\Da^{-\gamma}(v+\underline{\mathbf e})<\infty.
\end{equation}
We recall that for $a,b\in \O$, $\Da(a+b)\ge \Da(a)$. Let us prove that condition (\ref{eq:nessnotdirect}) implies that $\gamma>0$. Indeed, if this is not the case, that is if $\gamma\le 0$, then 
$$\Da^{\beta-\nu+\frac nr}(v)\Da^{-\gamma}(v+\underline{\mathbf e})\ge \Da^{\beta-\gamma-\nu+\frac nr}(v).$$
Hence for every $v\in \O$, we must have $\Da^{\beta-\gamma-\nu+\frac nr}(v)<\infty$, which is only possible if  $$\beta-\gamma-\nu+\frac nr=0.$$ But by (\ref{eq:cone21}), $\beta-\gamma-\nu+\frac nr=0$ only if $\mu=-q\alpha$. This is impossible since the latter used in (\ref{eq:cone23}) leads to $\frac nr-1<0$, which is not true. Hence $\gamma>0$.
\end{proof}
We finally prove the following.
\blem\label{lem:main3coneness}
Suppose $\nu\in \mathbb{R}$, and $1< p< \infty$. If the operator $S$ is bounded from $L^p_\nu(\O)$ into $L^\infty(\O)$, then the parameters satisfy the conditions (\ref{eq:cone31}), (\ref{eq:cone32}) and we have $\ga=\a+\beta+\frac nr-\frac \nu p$.
\elem
\begin{proof}
Again, the necessity of condition $\ga=\a+\beta+\frac nr-\frac \nu p$ is obtained as in the previous lemmas. 
To see that (\ref{eq:cone31}) and (\ref{eq:cone32}) holds, observe that the boundedness of $S$ 
from $L^p(\O, \Da^{\nu-\frac nr}(y)dy)$ to $L^\infty(\O)$ is equivalent to the boundedness of the adjoint $S^*$ of 
$S$ from $L^{1}(\O, \Da^{\nu-\frac nr}(y)dy)$ to $L^{p'}(\O, \Da^{\nu-\frac nr}(y)dy)$, where $S^*$ is given by
$$S^*g(v)=\Da^{\beta-\nu+\frac nr}(v)\int_{\O}\Da^{-\gamma}(y+v)g(y)\Da^{\alpha+\nu-\frac nr}(y)dy.$$
The proof  is then obtained as in the previous lemmas. 
\end{proof}
\subsection{Proof of Theorem \ref{thm:main1cone}, Theorem \ref{thm:main2cone} and Theorem \ref{thm:main3cone}}
Theorem \ref{thm:main1cone} follows from Lemma \ref{lem:conesuff1} and Lemma \ref{lem:main1coneness}. Theorem \ref{thm:main2cone} follows from Lemma \ref{lem:main2cone} and Lemma \ref{lem:main2coneness} while  Theorem \ref{thm:main3cone} is derived from Lemma \ref{lem:main3cone} and Lemma \ref{lem:main3coneness}.
\section{Off-diagonal estimates for the family $T_{\alpha,\beta,\gamma}^+$}
\subsection{The mixed norm case}
As said in the presentation of our results, following the proof of \cite[Theorem 1.1]{sehba}, one has that Theorem \ref{thm:main1}, Theorem \ref{thm:main2} and Theorem \ref{thm:main3} follow respectively from Theorem \ref{thm:main1cone}, Theorem \ref{thm:main2cone} and Theorem \ref{thm:main3cone}. 
\vskip .2cm
As a special case of Theorem \ref{thm:main1}, we have the following which extends the diagonal case result   ($q=s$) \cite[Theorem 4.3]{BBGNPR}  (see also \cite[Corollary 3.6]{sehba}).
\bcor\label{bounded}
Let $1<p<\infty$, $1<q\le s<\infty$, and assume that $\nu,\mu>\frac{n}{r}-1$. Then the operator $P_\nu^+$ is bounded from $L_\nu^{p,q}(\mathcal{D})$ to $L_\mu^{p,s}(\mathcal{D})$ if and only if $\frac{\nu}{q}=\frac{\mu}{s}$ and $1+\frac{\frac{n}{r}-1}{\mu}<q<1+\frac{\nu}{\frac{n}{r}-1}$.
\ecor 
If we restrict our problem to the boundedness of the positive Bergman projection $P_\beta^+$ from $L_\nu^{p,1}(\mathcal{D})$ to $L_\mu^{p,s}(\mathcal{D})$ with $s>1$, then we deduce the following from Theorem \ref{thm:main2}.
\bcor\label{cor:main2}
Let $1<p,s<\infty$,  and assume that $\beta,\mu>\frac{n}{r}-1$, and $\nu>0$. Then the following are equivalent.
\begin{itemize}
\item[(a)] The operator $P_\beta^+$ is bounded from $L_\nu^{p,1}(\mathcal{D})$ to $L_\mu^{p,s}(\mathcal{D})$
\item[(b)] The parameters satisfy $s\nu=\mu$ and $\frac{n}{r}-1<s(\beta-\nu)$.
\end{itemize}
\ecor
As observed in the second paragraph, the range of boundedness of the positive projection $P_\beta^+$ from $L_\nu^{p,q}(\mathcal{D})$ into $L_\mu^{p,s}(\mathcal{D})$ is not in general optimal for the boundedness of the projection $P_\beta$ (see for example \cite{BBPR,BBGNPR,sehba}). If in the corollary above, we choose $\beta$ large enough, we obtain that the optimal condition for the boundedness of $P_\beta^+$ is also optimal for $P_\beta$. In particular, putting $\beta=\nu+\frac{n}{r}-1$, we have the following full characterization for the Bergman projection $P_{\beta}$.
\bprop
Let $1<p<\infty$, $1<s<\infty$, and assume that $\mu>\frac{n}{r}-1$ and $\nu>0$. Then the following are equivalent.
\begin{itemize}
\item[(a)] The operator $P_{\nu+\frac{n}{r}-1}$ is bounded from $L_\nu^{p,1}(\mathcal{D})$ to $L_\mu^{p,s}(\mathcal{D})$ 
\item[(b)] The operator $P_{\nu+\frac{n}{r}-1}^+$ is bounded from $L_\nu^{p,1}(\mathcal{D})$ to $L_\mu^{p,s}(\mathcal{D})$
\item[(c)] The parameters satisfy $s\nu=\mu$.
\end{itemize}
\eprop 

\begin{proof}
We already know that the implication (b)$\Rightarrow$(a) holds. That (c)$\Rightarrow$(b) is a special case of Corollary \ref{cor:main2}. That (a)$\Rightarrow$(c) can be proved as in Lemma \ref{lem:main2coneness}.
\end{proof}

In the diagonal case, it is known that the Bergman projection $P_\gamma$ cannot be bounded on $L_\nu^{p,q}(\mathcal{D})$ for $q$ larger than $\tilde q_{\nu,p}$, where  $\tilde q_{\nu,p}=\frac{\nu +
\frac{n}{r}-1}{(\frac{n}{rp'}-1)_+} $ with $\tilde
q_{\nu,p}=\infty$ if $n/r\le p'$ (see \cite{BBPR,BBGR,BBGNPR,sehba}).  In the off-diagonal case, the same phenomenon is observed. Indeed we have the following.
\bprop Let $\mu,\nu,\gamma\in \mathbb{R}$, and $1\le p<\infty$ and $1\le q\le s<\infty$. If
$P_\gamma$ extends as a bounded operator from $L_\nu^{p,q}(\mathcal{D})$ into $L_\mu^{p,s}(\mathcal{D})$, then
$B_\mu(z,i\underline {\mathbf e})\in L_\mu^{p,s}$ and
$\Da^{\gamma-\nu}(\Im z)B_\mu(z,i\underline {\mathbf e})\in
L_\nu^{p',q'}$. The latter is equivalent to the following
conditions: $\mu>\frac{n}{r}-1$, $\left(\gamma+\frac{n}{r}\right)>\left(2\frac{n}{r}+1\right)\max\{\frac{1}{p},\frac{1}{p'}\}$, and $s>\frac{\mu+\frac{n}{r}-1}{\left(\gamma+\frac{n}{rp'}\right)_+}$
and 
$$\frac{\nu-\frac{n}{r}+1}{(\gamma-\frac{n}{r}+1)_+}<q<\tilde
q_{\nu,p}.$$ 
\eprop
\begin{proof} Let us denote by $P_\gamma^*$ the adjoint operator of $P_\gamma$ with respect to the
pairing $\langle\,,\rangle_\nu$. One easily checks that
 $$P_\gamma^* f(z)=\Da^{\gamma-\nu}(\Im
z)\int_{\mathcal{D}}B_\gamma(z,w)f(w)\Da^{\nu-n/r}(\Im w)dV(w),$$ 
$f\in L_{\mu+(\nu-\mu)s'-\frac{n}{r}}^{p',s'}(\mathcal{D}).$  Let $B_1(i\underline
{\mathbf e})$ be the Euclidean ball of radius 1 centered at
$i\underline{\mathbf e}$. Testing $P_\gamma$  with
$f_1(z)=\chi_{B_1(i\underline {\mathbf
e})}(z)\Da^{-\gamma+\frac{n}{r}}(\Im z)$ and $P_\mu^*$ with
$f_2(z)=\chi_{B_1(i\underline {\mathbf
e})}(z)\Da^{-\nu+\frac{n}{r}}(\Im z)$, we obtain with the help of the mean value property
that $P_\gamma f_1(z)=CB_\gamma(z,i\underline {\mathbf e})$ and $P_\gamma^*
f_2(z)=C\Da^{\gamma-\nu}(\Im z)B_\gamma(z,i\underline {\mathbf e})$.
It follows that we  should have  $B_\gamma(z,i\underline {\mathbf e})\in
L_\mu^{p,s}$ and $\Da^{\gamma-\nu}(\Im z)B_\gamma(z,i\underline {\mathbf
e})\in L_\nu^{p',q'}$. By Lemma \ref{lem:Apqfunction} this is equivalent
to $$\mu>\frac{n}{r}-1,\,\,\, \nu+(\gamma-\nu)q'>\frac{n}{r}-1,$$
$$\gamma+\frac{n}{r}>(2\frac{n}{r}-1)\max(\frac{1}{p'},\frac{1}{p})$$
and $$\gamma
+\frac{n}{r}>\max\{\frac{n}{rp'}+\frac{\nu+(\gamma-\nu)q'+\frac{n}{r}-1}{q'},\frac{n}{rp}+\frac{\mu+\frac{n}{r}-1}{s}\}.$$
That is, $\mu>\frac{n}{r}-1$,
$\gamma+\frac{n}{r}>(2\frac{n}{r}-1)\max(\frac{1}{p'},\frac{1}{p})$,
and $s>\frac{\mu+\frac{n}{r}-1}{\left(\gamma+\frac{n}{rp'}\right)_+}$, and
$\frac{\nu-\frac{n}{r}+1}{(\gamma-\frac{n}{r}+1)_+}<q<\tilde
q_{\nu,p}$.
\end{proof}
\subsection{Proof of Theorem \ref{thm:main4}}
We appeal again to Okikiolu result and follow the ideas in the proof of Theorem \ref{thm:main1cone}.
\begin{proof}[Proof of Theorem \ref{thm:main4}]
As said, we follow the proof of Theorem \ref{thm:main1cone}. Recall that 
$$\ga=\a+\beta+\frac nr-\frac 1p\left(\nu+\frac nr\right)+\frac 1q\left(\mu+\frac nr\right).$$
We put
$$\omega=\a+\beta-\ga-\nu=-\left[\frac{\nu+\frac nr}{p'}+\frac{\mu+\frac nr}{q}\right]<0.$$
The inequality (\ref{eq:main41}) is equivalent to $\frac{\nu-\frac nr+1}{p}+\frac{1}{q}(\frac nr-1)<\beta+1.$ Hence,
 $$\beta-\nu+\frac nr+\frac 1{p'}(\nu-\frac nr+1)-\frac{1}{q}(\frac nr-1)>0.$$ Multiplying this last 
inequality by $\omega<0$ yields $$(\beta-\nu+\frac nr)\omega+\frac 1{p'}(\nu-\frac nr+1)\omega-\frac{\omega}{q}(\frac nr-1)<0,$$ i.e.
\Bea\label{4.3}
\frac{\nu-\frac nr+1}{p'}\omega-\frac{\beta-\nu+\frac nr}{p'}(\nu+\frac nr)<\frac{\beta-\nu+\frac nr}{q}(\mu+\frac nr)+\frac{\omega}{q}(\frac nr-1).
\Eea

From the inequality (\ref{eq:main42}), we have $\frac{\mu-\frac nr+1}{q}>-\a+\frac{1}{p'}(\frac nr-1)$. Multiplying this last 
inequality by $\omega$ yields $\frac{\mu-\frac nr+1}{q}\omega<-\a\omega+\frac{\omega}{p'}(\frac nr-1)$, i.e.
\Bea\label{4.4}
\frac{\mu-\frac nr+1}{q}\omega-\a\frac{\mu+\frac nr}{q}<\a\frac{\nu+\frac nr}{p'}+\frac{\omega}{p'}(\frac nr-1).
\Eea
The inequalities (\ref{4.3}) and (\ref{4.4}) yield the existence of two real numbers $u$ and $v$ such that
\Bea\label{u4}
\left\{\Ba{l}\frac{\nu-\frac nr+1}{p'}\omega-\frac{\beta-\nu+\frac nr}{p'}(\nu+\frac nr)<u\omega+(\beta-\nu+\frac nr)(v-u)\\
\quad\quad\quad\quad\quad\quad\quad\quad\quad\quad\quad\quad\quad\quad\quad\quad<\frac{\beta-\nu+\frac nr}{q}(\mu+\frac nr)+\frac{\omega}{q}(\frac nr-1)\\\\
\frac{\mu-\frac nr+1}{q}\omega-\a\frac{\mu+\frac nr}{q}<v\omega+\a(u-v)<\a\frac{\nu+\frac nr}{p'}+\frac{\omega}{p'}(\frac nr-1).\Ea\right.
\Eea

Now, (\ref{u4}) is equivalent to
\Bea\label{u41}
\left\{\Ba{l}-\frac{\beta-\nu+\frac nr}{\omega}\left[-\frac{\mu+\frac nr}{q}-u+v\right]+\frac{1}{q}(\frac nr-1)<u<\frac{\nu-\frac nr+1}{p'}+\frac{\beta-\nu+\frac nr}{\omega}\left[-\frac{\nu+\frac nr}{p'}+u-v\right]\\\\
-\frac{\a}{\omega}\left[-\frac{\nu+\frac nr}{p'}+u-v\right]+\frac{1}{p'}(\frac nr-1)<v<\frac{\mu-\frac nr+1}{q}+\frac{\a}{\omega}\left[-\frac{\mu+\frac nr}{q}-u+v\right].\Ea\right.
\Eea
Let
$$t=\frac{-\frac{\nu+\frac nr}{p'}+u-v}{\omega}$$
then
$$1-t=\frac{-\frac{\mu+\frac nr}{q}-u+v}{w}.$$
Since $\omega<0,$ we choose $u$ and $v$ such that $0<v-u<\frac{\mu+\frac nr}{q}.$ Thus, we have $0<t<1.$
Moreover, (\ref{u41}) becomes
\Bea\label{u42}
\left\{\Ba{l}-(\beta-\nu+\frac nr)(1-t)+\frac{1}{q}(\frac nr-1)<u<\frac{\nu-\frac nr+1}{p'}+(\beta-\nu+\frac nr)t\\\\
-\a t+\frac{1}{p'}(\frac nr-1)<v<\frac{\mu-\frac nr+1}{q}+\a(1-t).\Ea\right.
\Eea

\vskip .3cm
We now observe that the positive kernel of the operator $T^+$ with respect to $\Da^{\nu-\frac nr}(\Im w)dV(w)$ is given by
$$K(z,w)=\Da^\a(\Im z)\left|\Da^{-\ga-\frac nr}\left(\frac{z-\overline{w}}{i}\right)\right|\Da^{\beta-\nu+\frac nr}(\Im w).$$
Let us consider the following positive functions $\psi_1(w)=\Da^{-u}(\Im w)$ and $\psi_2(z)=\Da^{-v}(\Im z).$ Then
\Beas
J_1&=&\int_{\mathcal{D}} K(z,w)^{tp'}\psi_1(w)^{p'}dV_\nu(w)\\
&=&\Da^{tp'\a}(\Im z)\int_{\mathcal{D}}\left|\Da^{-tp'\left(\ga+\frac nr\right)}\left(\frac{z-\overline{w}}{i}\right)\right|\Da^{tp'(\beta-\nu+\frac nr)-p'u+\nu-\frac nr}(\Im w)dV(w).
\Eeas
The last integral above converges because from the right inequality in (\ref{u42}) involving $u$,
we have $tp'(\beta-\nu+\frac nr)-p'u+\nu>\frac nr-1$ and also
\Beas 
I &:=& -tp'(\ga+\frac nr)+tp'(\beta-\nu+\frac nr)-p'u+\nu\\ &=& -tp'(\ga-\beta+\nu)-p'u+\nu\\
&=& -tp'(\a+\nu+\frac nr-\frac{1}{p}(\nu+\frac nr)+\frac{1}{q}(\mu+\frac nr))-p'u+\nu\\ &=& -tp'(\a-\omega)-p'u+\nu\\
&=& -\a tp'+t\omega p'-p'u+\nu\\ &=& -\a tp'+p'(-\frac{1}{p'}(\nu+\frac nr)+u-v)-p'u+\nu\\
&=& -\a tp'-p'v-\frac nr<-2\frac nr+1
\Eeas
where the last inequality above follows from the left inequality in (\ref{u42}) involving $v.$
 It follows using Lemma \ref{lem:Apqfunction} that
\Beas
J_1&=&\Da^{tp'\a}(\Im z)\int_{\mathcal{D}}\left|\Da^{-tp'\left(\ga+\frac nr\right)}\left(\frac{z-\overline{w}}{i}\right)\right|\Da^{tp'(\beta-\nu+\frac nr)-p'u+\nu-\frac nr}(\Im w)dV(w)\\
&=&C_1\Da^{tp'\a-tp'\left(\ga+\frac nr\right)+tp'(\beta-\nu+\frac nr)-p'u+\nu+\frac nr}(\Im z)\\
&=&C_1\Da^{-p'v}(\Im z)=C_1\phi_2(z)^{p'}.
\Eeas

Next,
\Beas
J_2&=&\int_{\mathcal{D}} K(z,w)^{(1-t)q}\psi_2(z)^{q}dV_\mu(z)\\
&=&\Da^{(1-t)q(\beta-\nu+\frac nr)}(\Im w)\int_{\mathcal{D}}\left|\Da^{-(1-t)q(\ga+\frac nr)}\left(\frac{z-\overline{w}}{i}\right)\right|\Da^{(1-t)q\a-qv+\mu-\frac nr}(\Im z)dV(z).
\Eeas
The last integral above converges since from the right inequality in (\ref{u42}) involving $v$ ,
we have $(1-t)q\a-qv+\mu>\frac nr-1$ and also , we have
\Beas
I &:=&-(1-t)q(\ga+\frac nr)+(1-t)q\a-qv+\mu\\ &=& -(1-t)q(\ga-\a+\frac nr)-qv+\mu\\
&=& -(1-t)q\left[\beta+2\frac nr-\frac{1}{p}(\nu+\frac nr)+\frac{1}{q}(\mu+\frac nr)\right]-qv+\mu\\
&=& -(1-t)q\left[\beta-\nu+\frac nr+\frac{1}{p'}(\nu+\frac nr)+\frac{1}{q}(\mu+\frac nr)\right]-qv+\mu\\ &=& -(1-t)q\left[\beta-\nu+\frac nr-\omega\right]-qv+\mu\\ &=& (1-t)q\omega-(1-t)q(\beta-\nu+\frac nr)-qv+\mu\\ &=& q(-\frac{1}{q}(\mu+\frac nr)-u+v) -(1-t)q(\beta-\nu+\frac nr)-qv+\mu\\ &=& -qu-(1-t)q(\beta-\nu+\frac nr)-\frac nr<-2\frac nr+1.
\Eeas
where the last inequality above follows from the left inequality in (\ref{u42}) involving $u.$
It follows from Lemma \ref{lem:Apqfunction} that
\Beas
J_2&=&\Da^{(1-t)q(\beta-\nu+\frac nr)}(\Im w)\int_{\mathcal{D}}\left|\Da^{-(1-t)q(\ga+\frac nr)}\left(\frac{z-\overline{w}}{i}\right)\right|\Da^{(1-t)q\a-qv+\mu-\frac nr}(\Im z)dV(z)\\
&=&C_2\Da^{(1-t)q(\beta-\nu+\frac nr)-(1-t)q(\ga+\frac nr)+(1-t)q\a-qv+\mu+\frac nr}(\Im w)\\
&=&C_2\Da^{-qu}(\Im w)=C_2\phi_1(w)^{q}
\Eeas
Thus by the Okikiolu test, we conclude that $T^+:L^p_\nu(\mathcal{D})\to L^q_\mu(\mathcal{D})$ is bounded.
\end{proof}

\bibliographystyle{plain}

\end{document}